**Some remarks on the power product expansion of the q-exponential series**

*Johann Cigler*

Fakultät für Mathematik, Universität Wien

**Abstract.** We give an overview about the power product expansion of the exponential series and derive some $q-$analogs.

## 1. Introduction

Each formal power series $f(x)=\sum_{n=0}^{\infty} a_n x^n$ with $a_0=1$ has a representation as an infinite product of the form

$$f(x)=\prod_{n=1}^{\infty}\left(1+g_n x^n\right), \tag{1}$$

a so-called power product expansion.

More precisely there are uniquely determined coefficients $g_j$ such that $\prod_{k=1}^{n}\left(1+g_k x^k\right)=\sum_k b_{n,k} x^k$ with $b_{n,k}=a_k$ for $k\leq n.$ To see this observe that

$$\prod_{k=1}^{n}\left(1+g_k x^k\right)=\left(1+g_n x^n\right)\sum_k b_{n-1,k} x^k=\sum_{k=0}^{n-1} a_k x^k+\left(b_{n-1,n}+g_n\right)x^n+\cdots.$$

Therefore, there is a uniquely determined $g_n$ such that $b_{n-1,n}+g_n=a_n.$

Such expansions have been studied in a number of papers, cf. [4] and [5] and the references cited there.

A simple example is

$$\sum_{n\geq 0} x^n=\prod_{n\geq 0}\left(1+x^{2^n}\right). \tag{2}$$

Here both sides are convergent for $|x|<1.$ Since the zeros of the right-hand side are not in the domain of convergence there is no contradiction to the fact that the left-hand side has no zeros.

In this note I consider the case of the exponential series

$$\exp(x)=\sum_{n\geq 0}\frac{x^n}{n!}=\prod_{n\geq 1}\left(1+e_n x^n\right) \tag{3}$$

and its $q-$analogs in some detail.

My interest in this problem has been aroused by the article “A dream of a (number-) sequence” by Gottfried Helms [6].

## 2. Some background information

Taking logarithms of (3) yields

$$x=\sum_{k\geq 1}\log\left(1+e_k x^k\right)=\sum_{k\geq 1}\sum_{d\geq 1}(-1)^{d-1}\frac{1}{d}e_k^d x^{kd}=\sum_{n\geq 1}x^n\sum_{d|n}(-1)^{d-1}\frac{e_{n/d}^d}{d}. \tag{4}$$

Therefore, we get for $n>1$

$\sum_{d|n}(-1)^{d-1}\frac{e_{n/d}^d}{d}=0$ which gives

$$e_n=\sum_{d|n,d>1}(-1)^d\frac{e_{n/d}^d}{d}. \tag{5}$$

The first terms are

$$\left(e_n\right)_{n\geq 1}=\left(1,\frac{1}{2},-\frac{1}{3},\frac{3}{8},-\frac{1}{5},\frac{13}{72},-\frac{1}{7},\frac{27}{128},\cdots\right) \tag{6}$$

O. Kolberg [7] used this formula to show that $e_p=-\frac{1}{p}$ for a prime number $p\geq 3$ and that $(-1)^n e_n>0$ for $n>1$.

The first assertion follows immediately from (5). For the second assertion we consider $\exp(-x)=\prod_{n\geq 1}\left(1+a_n x^n\right)$ and show that $0<a_n<\frac{2}{n}$ for $n>1$.

Let me reproduce his argument. The inequality is true for $n<20$ by direct computation. Let now $n\geq 20$ and suppose that $0<a_m<\frac{2}{m}$ holds for $m<n$. Then

$$\left|a_n-\frac{1}{n}\right|<\left|\sum_{d|n,1<d<n}(-1)^d\frac{1}{d}a_{n/d}^d\right|<\frac{1}{2}\left(\frac{4}{n}\right)^2+\frac{1}{3}\left(\frac{6}{n}\right)^3+\frac{1}{4}\left(\frac{8}{n}\right)^4+\sum_{d|n,5\leq d\leq\frac{n}{4}}\frac{1}{d}a_{n/d}^5+\frac{2}{n}a_2^{\frac{n}{2}}+\frac{3}{n}a_3^{\frac{n}{3}}.$$

For $d\leq\frac{n}{4}$ we have $k=\frac{n}{d}\geq 4$ and $a_k<\frac{2}{k}$. This implies $\sum_{d|n,5\leq d\leq\frac{n}{4}}\frac{1}{d}a_{n/d}^5<\sum_{k\geq 4}\frac{k}{n}\left(\frac{2}{k}\right)^5$

For $n\geq 20$ we get $n\left(\frac{1}{2}\left(\frac{4}{n}\right)^2+\frac{1}{3}\left(\frac{6}{n}\right)^3+\frac{1}{4}\left(\frac{8}{n}\right)^4+\frac{2}{n}a_2^{\frac{n}{2}}+\frac{3}{n}a_3^{\frac{n}{3}}\right)<0.72.$

This implies $\left|a_n-\frac{1}{n}\right|<\frac{1}{n}\left(0.72+\sum_{k\geq 4}k\left(\frac{2}{k}\right)^5\right)<\frac{1}{n}\left(0.72+0.24\right)<\frac{1}{n}.$

It also follows that the expansion

$$\frac{\exp(-x)}{1-x} = \prod_{n\geq 2}\left(1+a_n x^n\right) \tag{7}$$

converges for $|x|<1$.

L. Carlitz [2] showed that $e_n = \frac{c_n}{n!}$ with integers $c_n$ and derived some arithmetic properties of $c_n$, for example that $c_n \equiv 0 \bmod p$ if $n > p$ is relatively prime to the prime number $p$. We shall give a different proof in Corollary 4.4.

The first terms of the sequence $(c_n)$ are

$$\left(c_n\right)_{n\geq 1} = \left(1,1,-2,9,-24,130,-720,8505,\cdots\right). \tag{8}$$

**Remark**

The sequence $\left(0,1,2,9,24,130,720,8505,\cdots\right)$ of coefficients of (7) also occurs as the dimensions of representations by Witt vectors (cf. [1] and OEIS [8], A006973).

Since $a_p = \frac{1}{p}$ for a prime number $p$ we get $c_p = -(p-1)!$.

Moreover J. Borwein and S. Lou [1] proved that $|c_n| \leq (n-1)!$ for odd $n$, for example $|c_9| = 35840 < 8! = 40320$, and that $c_n \geq (n-1)!$ for even $n$, for example $c_4 = 9 > 3! = 6$ or $c_6 = 130 > 5! = 120$.

A look at OEIS [8] suggests the sequence OEIS [8], A067911= $\left(u_n = \prod_{k=1}^{n}\gcd(k,n) = \prod_{d|n} d^{\varphi\left(\frac{n}{d}\right)}\right)_{n\geq 1}$ as another choice for the numerators $1,2,3,8,5,72,7,128,\cdots$.

Putting $e_n = \frac{r_n}{u_n}$ in (5) we get

$$r_n = \sum_{d|n,d>1} (-1)^d \frac{u_n}{d u_{n/d}^d} r_{n/d}^d. \tag{9}$$

Now we see by induction that all $r_n$ are integers if we can show that $\frac{u_n}{d u_{n/d}^d}$ is an integer for each divisor $d$ of $n$.

In order to show this we consider the $d$ intervals $j\frac{n}{d}+1,\cdots,(j+1)\frac{n}{d}$, $0\le j<d$, where in the last interval the number $n$ is replaced by $\frac{n}{d}$. Since for each $i$ $\gcd(i,n/d)\big|\gcd(i,n)$ we see that

$u_{n/d}^{d}\left|\frac{u_n}{d}\right.$.

The first terms are

$$\left(r_n\right)_{n\ge 1}=\left(1,1,-1,3,-1,13,-1,27,\cdots\right). \tag{10}$$

Note that $r_n$ and $u_n$ need not be relatively prime. For example, $\gcd\left(u_{12},r_{12}\right)=3$.

Let us mention some special cases of (9): For a prime $p>2$ we get $r_p=-1$ and $r_{p^2}=1-p^{p-1}$ and for a product of two different primes $n=pq$ we get $u_{pq}=p^q q^p$ and

$$\begin{aligned} r_{pq}&=\sum_{\substack{d|n\\ d>1}}(-1)^d\frac{u_n}{du_{n/d}^d}r_{n/d}^d=-\frac{u_{pq}}{pu_q^p}r_q-\frac{u_{pq}}{qu_p^q}r_p-\frac{u_{pq}}{pq}=u_{pq}\left(\frac{1}{pq^p}+\frac{1}{qp^q}-\frac{1}{pq}\right)\\ &=p^q q^p\left(\frac{1}{pq^p}+\frac{1}{qp^q}-\frac{1}{pq}\right)=p^{q-1}+q^{p-1}-p^{q-1}q^{p-1}. \end{aligned}$$

### 3. Connections with Pascal's triangle.

**3.1.** Let $P_n=\left(\binom{i}{j}\right)_{i,j=0}^{n-1}$ be the $n\times n-$Pascal matrix and let $H_n=\left(h(i,j)\right)_{i.j=0}^{n-1}$ with $h(i,i-1)=i$ and $h(i,j)=0$ else, i.e. $h(i,j)=i\left[i-j=1\right]$ by using Iverson's convention: [P]=1 if property P is true and [P]=0 else. Then $\frac{H_n^k}{k!}=H_{n,k}=\left(\binom{i}{k}\left[i-j=k\right]\right)_{i,j=0}^{n-1}$.

Since $H_n^k=0$ for $k\ge n$ we get

$$P_n=\sum_{k=0}^{n-1}H_{n,k}=\sum_{k=0}^{n-1}\frac{1}{k!}H_n^k=\sum_{k\ge 0}\frac{1}{k!}H_n^k=\exp\left(H_n\right). \tag{11}$$

From (3) we get

$$\begin{aligned} P_n&=\left(I_n+\frac{c_1}{1!}H_n^1\right)\left(I_n+\frac{c_2}{2!}H_n^2\right)\cdots\left(I_n+\frac{c_{n-1}}{(n-1)!}H_n^{n-1}\right)\\ &=\left(I_n+c_1H_{n,1}\right)\left(I_n+c_2H_{n,2}\right)\cdots\left(I_n+c_{n-1}H_{n,n-1}\right). \end{aligned} \tag{12}$$

For example

$$P_4=\left(I_4+H_{4,1}\right)\left(I_4+H_{4,2}\right)\left(I_4-2H_{4,3}\right)=\begin{pmatrix}1&0&0&0\\1&1&0&0\\1&2&1&0\\1&3&3&1\end{pmatrix}=\begin{pmatrix}1&0&0&0\\1&1&0&0\\0&2&1&0\\0&0&3&1\end{pmatrix}\begin{pmatrix}1&0&0&0\\0&1&0&0\\1&0&1&0\\0&3&0&1\end{pmatrix}\begin{pmatrix}1&0&0&0\\0&1&0&0\\0&0&1&0\\-2&0&0&1\end{pmatrix}.$$

From this representation we see again that the numbers $c_n$ are integers.

If we set $\left(I_n+c_1H_{n,1}\right)\left(I_n+c_2H_{n,2}\right)\cdots\left(I_n+c_kH_{n,k}\right)=\left(g_k(i,j)\right)_{i,j=0}^{n-1}$ then $g_k(i,0)=1$ for $i\le k.$

If it is already proved that $c_i$ is an integer for $1\le i\le k-1,$ then from $\left(g_k(i,j)\right)_{i,j=0}^{n-1}=\left(g_{k-1}(i,j)\right)_{i,j=0}^{n-1}\left(I_n+c_kH_{n,k}\right)$ we see that $g_k(k,0)=g_{k-1}(k,0)+c_k$. Choosing $c_k=1-g_{k-1}(k,0)$ we get $g_k(k,0)=1$ and all $g_k(i,j)$ are integers.

**3.2.** For later applications let us state a slight generalization. Let $m$ be a positive integer and consider the matrices $H_{n,k}^{(m)}=\left(h_k^{(m)}(i,j)\right)_{i,j=0}^{n-1}=\left(\binom{\left\lfloor \frac{i}{m}\right\rfloor}{k}[i-j=mk]\right)_{i,j=0}^{n-1}$ whose entries $h_k^{(m)}(i,j)$ satisfy

$h_k^{(m)}(i,i-mk)=\binom{\left\lfloor \frac{i}{m}\right\rfloor}{k}$ and $h_k^{(m)}(i,j)=0$ else.

They satisfy $H_{n,k}^{(m)}=kH_{n,k-1}^{(m)}H_{n,1}^{(m)}$ because

$$\sum_\ell h_{k-1}^{(m)}(i,\ell)h_1^{(m)}(\ell,i-mk)=h_{k-1}^{(m)}(i,i-m(k-1))h_1^{(m)}(i-m(k-1),i-mk)$$

$$=\binom{\left\lfloor \frac{i}{m}\right\rfloor}{k-1}\binom{\left\lfloor \frac{i}{m}\right\rfloor-k+1}{1}=k\binom{\left\lfloor \frac{i}{m}\right\rfloor}{k}$$

and $\sum_\ell h_{k-1}^{(m)}(i,\ell)h_1^{(m)}(\ell,j)=0$ else. Therefore we get $H_{n,k}^{(m)}=\dfrac{\left(H_n^{(m)}\right)^k}{k!}$ for $H_n=H_{n,1}$

and

$$P_n^{(m)}=\sum_{k\ge0}H_{n,k}^{(m)}=\exp\left(H_n^{(m)}\right)=\prod_{k\ge1}\left(I_n+c_kH_{n,k}^{(m)}\right). \tag{13}$$

For $m=2$ we get the " doubled" Pascal triangle (OEIS [8], A178112).

## 4. Power product expansion of the $q-$exponential series

Let us now consider the $q-$exponential series

$$\exp_q(x)=\sum_{n\geq 0}\frac{x^n}{[n]!} \tag{14}$$

and its counterpart

$$\operatorname{Exp}_q(x)=\sum_{n\geq 0}q^{\binom{n}{2}}\frac{x^n}{[n]!}. \tag{15}$$

As usual we set $[n]=[n]_q=1+q+\cdots+q^{n-1}$ and $[n]!=[1][2]\cdots[n]$.

Depending on the context $q$ is either a complex number or an indeterminate. The needed results about $q-$series may for example be found in [3].

**Theorem 4.1**

*The coefficients* $e_n(q)$ *of the power product expansion of*

$$\exp_q(x)=\prod_{n\geq 1}\left(1+e_n(q)x^n\right) \tag{16}$$

*are given by*

$$e_n(q)=\sum_{d|n,d>1}(-1)^d\frac{e_{n/d}(q)^d}{d}+\frac{(1-q)^{n-1}}{n[n]}. \tag{17}$$

*The coefficients* $E_n(q)$ *of the expansion*

$$\operatorname{Exp}_q(x)=\exp_{q^{-1}}(x)=\frac{1}{\exp_q(-x)}=\prod_{n\geq 1}\left(1+E_n(q)x^n\right) \tag{18}$$

*are given by*

$$E_n(q)=\sum_{d|n,d>1}(-1)^d\frac{E_{n/d}(q)^d}{d}+\frac{(q-1)^{n-1}}{n[n]}. \tag{19}$$

*They satisfy*

$$e_{2n+1}(q)=E_{2n+1}(q)=e_{2n+1}\left(\frac{1}{q}\right). \tag{20}$$

**Proof**

For $0\leq q<1$ we have $\exp_q(x)=\prod_{k=0}^{\infty}\left(1-q^k(1-q)x\right)^{-1}$ and therefore

$$\log \exp_q(x) = -\sum_{k=0}^{\infty} \log\left(1-q^k(1-q)x\right) = \sum_{k=0}^{\infty}\sum_{n=1}^{\infty} \frac{\left(q^k(1-q)x\right)^n}{n}$$
$$= \sum_{n=1}^{\infty} \frac{\left((1-q)x\right)^n}{n} \frac{1}{1-q^n} = \sum_{n=1}^{\infty} \frac{(1-q)^{n-1}}{n[n]} x^n. \tag{21}$$

This gives (17). In the same way we get (19). Comparing both formulas gives (20).

The first terms of the sequence $\left(e_n(q)\right)$ are

$$e_1(q)=1,\ e_2(q)=\frac{1}{1+q},\ e_3(q)=-\frac{q}{[3]},\ e_4(q)=\frac{\left(1+q^2+q^3\right)}{[2][4]},\ e_5(q)=-\frac{q\left(1-q+q^2\right)}{[5]},$$
$$e_6(q)=\frac{q^2\left(1+3q+2q^2+2q^3+2q^4+2q^5+q^6\right)}{[2]^2[3][6]},\ e_7(q)=-\frac{q\left(1-q+q^2\right)^2}{[7]},\cdots.$$

The first terms of the sequence $\left(E_n(q)\right)$ are

$$E_1(q)=1,\ E_2(q)=\frac{q}{1+q},\ E_3(q)=-\frac{q}{[3]},\ E_4(q)=\frac{q\left(1+q+q^3\right)}{[2][4]},\ E_5(q)=-\frac{q\left(1-q+q^2\right)}{[5]},$$
$$E_6(q)=\frac{q\left(1+2q+2q^2+2q^3+2q^4+3q^5+q^6\right)}{[2]^2[3][6]},\ E_7(q)=-\frac{q\left(1-q+q^2\right)^2}{[7]},\cdots.$$

Let now $u_n(q)=\prod_{j=1}^{n}\gcd\left([j],[n]\right)$ be the product of the polynomial greatest divisors of the polynomials $[j]=1+q+\cdots+q^{j-1}$ and $[n]=1+q+\cdots+q^{n-1}$. Since $\left(q^d-1\right)\Big|\left(q^n-1\right)$ if and only if $d\,|\,n$ we see that $\lim_{q\to 1} u_n(q)=u_n$.

Let us write $e_n(q)=\frac{r_n(q)}{u_n(q)}$. Then we get as above

$$r_n(q)=\sum_{d|n,d>1}(-1)^d \frac{u_n(q)}{du_{n/d}^d(q)} r_{n/d}^d(q)+\frac{(1-q)^{n-1}u_n(q)}{n[n]}. \tag{22}$$

The first terms of $\left(r_n(q)\right)$ are

$$1,1,-q,1+q^2+q^3,-q\left(1-q+q^2\right),q^2\left(1+3q+2q^2+2q^3+2q^4+2q^5+q^6\right),-q\left(1-q+q^2\right)^2,\cdots.$$

**Remark**

For $q=0$ the series $\exp_q(x)$ reduces to $1+x+x^2+\cdots=\dfrac{1}{1-x}$ and therefore (16) reduces to (2).

Comparing coefficients we get $r_{2^n}(0)=1$ and $r_n(0)=0$ else.

**Theorem 4.2**

*The polynomials* $(-1)^n r_n(q)$ *have integer coefficients and leading coefficient* $1$ *for* $n>1$.

**Proof**

Let us first show that $(-1)^n r_n(q)$ has leading coefficient $1$ for $n>1$.

The highest terms of $q$ in $(-1)^n r_n(q)$ occur in

$$\frac{u_n(q)}{n}-\frac{(q-1)^{n-1}u_n(q)}{n[n]}=\frac{u_n(q)}{n}\left(1-\frac{(q-1)^{n-1}}{[n]}\right)=\frac{u_n(q)}{[n]}\frac{[n]-(q-1)^{n-1}}{n}.$$

The coefficient of the leading term of $\dfrac{u_n(q)}{[n]}$ is $1$ and the leading term of $\dfrac{[n]-(q-1)^{n-1}}{n}$ is $q^{n-2}$.

Next we show that $r_n(q)$ is a polynomial in $q$ with integer coefficients.

We show first that $e_n(q)=\dfrac{c_n(q)}{[n]!}$ where $c_n(q)$ is a polynomial in $q$ with integer coefficients.

Let $P_n(q)=\left(\begin{bmatrix} i \\ j \end{bmatrix}_q\right)_{i.j=0}^{n-1}$ and let $H_n(q)=\left(\begin{bmatrix} i \\ 1 \end{bmatrix}_q [i-j=1]\right)_{i,j=0}^{n-1}$.

Then $H_n^k(q)=\left([k]!\begin{bmatrix} i \\ k \end{bmatrix}_q [i-j=k]\right)_{i,j=0}^{n-1}$.

This implies the well-known $q-$analog of (11)

$$\exp_q\left(H_n(q)\right)=\sum_{k=0}^{\infty}\frac{H_n^k(q)}{[k]!}=\left(\begin{bmatrix} i \\ j \end{bmatrix}\right)_{i,j=0}^{n-1}=P_n(q). \tag{23}$$

If we write $\dfrac{H_n^k(q)}{[k]!}=H_{n,k}(q)$ then

$$P_n(q)=\prod_{k\geq 1}\left(I_n+e_k(q)H_n^k(q)\right)=\prod_{k\geq 1}\left(I_n+c_k(q)H_{n,k}(q)\right). \tag{24}$$

We now show that each $c_n(q)$ is a polynomial in $q$ with integer coefficients. Since $c_1(q)=1$ this is true for $k=1$. Let $\prod_{j=0}^{k}\left(I_n+c_j(q)H_{n,j}(q)\right)=\left(g_{k,i,j}(q)\right)_{i,j=0}^{n-1}$.

Assume that all $g_{k-1,i,j}(q)$ are polynomials with integer coefficients. We know that $g_{k,0,k}=1$ since for $n=k$ we get $P_k(q)$. Then $1=g_{k,0,k}=g_{k-1,0,k}+g_{k-1,k,k}c_k(q)=g_{k-1,0,k}+c_k(q)$ shows that $c_k(q)$ has integer coefficients and therefore that all $g_{k,i,j}(q)$ are polynomials in $q$ with integer coefficients.

Since $\frac{1}{[k]!}$ is a formal power series with integer coefficients we see that

$r_k(q)=\frac{u_k(q)c_k(q)}{[k]!}\in\mathbb{Z}[q]$ is also a polynomial with integer coefficients.

The first terms of the sequence $\left(r_n(q)\right)_{n\geq 1}$ are $1, 1, -q, 1+q^2+q^3, -q\left(1-q+q^2\right)$, $q^2\left(1+3q+2q^2+2q^3+2q^4+2q^5+q^6\right), -q\left(1-q+q^2\right)^2,\cdots$.

For a prime number $n=p$ we get $u_p(q)=[p]$ and

$$r_p(q)=-\frac{u_p(q)}{p}+\frac{(1-q)^{p-1}u_p(q)}{p[p]}=-\frac{[p]}{p}+\frac{(1-q)^{p-1}}{p}=\frac{1}{p(1-q)}\left(-1+q^p+(1-q)^p\right)$$

$$=\frac{\sum_{j=1}^{p-1}\binom{p}{j}(-q)^j}{p(1-q)}=\frac{\sum_{j=1}^{\frac{p-1}{2}}\binom{p}{j}(-q)^j\left(1-q^{p-2j}\right)}{p(1-q)}\in\mathbb{Z}[q].$$

For the sequence $\left(c_n(q)\right)_{n\geq 1}$ we get more information.

**Theorem 4.3**

*Let $n\geq m$ be positive integers and $\zeta_m=e^{\frac{2\pi i}{m}}$. If $n=mk$ we get $c_{mk}\left(\zeta_m\right)=c_k$, if $n$ is not a multiple of $m$ then $c_n\left(\zeta_m\right)=0$.*

**Proof**

Since $\frac{[km+r]_q!}{[km]_q!}=\frac{1-q^{km+1}}{1-q}\frac{1-q^{km+2}}{1-q}\cdots\frac{1-q^{km+r}}{1-q}$ for $q\to\zeta_m$ reduces to $\frac{1-\zeta_m}{1-\zeta_m}\frac{1-\zeta_m^{\,2}}{1-\zeta_m}\cdots\frac{1-\zeta_m^{\,r}}{1-\zeta_m}=[r]_{\zeta_m}!$ we see that for $q\to\zeta_m$

$$H_{n,km+r}\left(\zeta_m\right)=\lim_{q\to\zeta_m}\frac{H_n^{km+r}(q)}{[km+r]!}=\lim_{q\to\zeta_m}\frac{H_n^{km}(q)}{[km]!}\frac{[km]!}{[km+r]!}H_n^r(q)=H_{n,km}\left(\zeta_m\right)\frac{H_n^r(\zeta_m)}{[r]_{\zeta_m}!}.$$

Therefore we get

$$P_n(\zeta_m) = \sum_{j=0}^{m-1} \frac{H_n^j(\zeta_m)}{[j]_{\zeta_m}!} \sum_{k\geq 0} H_{n,km}(\zeta_m). \tag{25}$$

Since

$$\begin{bmatrix} nm+r \\ km \end{bmatrix} = \frac{[nm+r]!}{[km]![(n-k)m+r]!} = \frac{[nm]!}{[km]![(n-k)m]!} \frac{[nm+r]}{[(n-k)m+r]} \frac{[nm+r-1]}{[(n-k)m+r-1]} \cdots \frac{[nm+1]}{[(n-k)m+1]}$$

and $\dfrac{[nm+i]}{[(n-k)m+i]} = \dfrac{1-q^{nm+i}}{1-q^{(n-k)m+i}} = 1$ for $q = \zeta_m$ and $0 < i < m$ we see that

$$\lim_{q\to\zeta_m} \begin{bmatrix} nm+r \\ km \end{bmatrix}_q = \lim_{q\to\zeta_m} \begin{bmatrix} nm \\ km \end{bmatrix}_q = \lim_{q\to\zeta_m} \frac{1-q^{mn}}{1-q^{km}} \frac{1-q^{mn-1}}{1-q^{km-1}} \cdots \frac{1-q^{m(n-k)+1}}{1-q^1} = \lim_{q\to\zeta_m} \begin{bmatrix} n \\ k \end{bmatrix}_{q^m} = \lim_{q\to 1} \begin{bmatrix} n \\ k \end{bmatrix}_q = \binom{n}{k}.$$

Therefore we get $\begin{bmatrix} i \\ km \end{bmatrix}_{q=\zeta_m} = \dbinom{\left\lfloor \frac{i}{m} \right\rfloor}{k}$ and

$$H_{n,km}(\zeta_m) = \left( \begin{bmatrix} i \\ km \end{bmatrix}_{\zeta_m} [i-j=km] \right)_{i,j=0}^{n-1} = \left( \binom{\left\lfloor \frac{i}{m} \right\rfloor}{k} [i-j=km] \right)_{i,j=0}^{n-1} = H_{n,k}^{(m)} = \frac{\left(H_n^{(m)}\right)^k}{k!}.$$

Thus (25) gives

$$\left( \sum_{j=0}^{m-1} \frac{H_n^j(\zeta_m)}{[j]_{\zeta_m}!} \right)^{-1} P_n(\zeta_m) = \exp\left(H_n^{(m)}\right) = \prod_{k\geq 1} \left( I_n + c_k H_{n,km}(\zeta_m) \right). \tag{26}$$

On the other hand we know that

$$\prod_{j=1}^{m-1} \left( I_n + c_j(q) H_{n,j}(q) \right) = \sum_{j=0}^{m-1} \frac{H_n^j(q)}{[j]!} + \sum_{j\geq m} b_j(q) H_n^j(q) \tag{27}$$

for some polynomials $b_j(q) \in \mathbb{Z}[q]$.

For $q = \zeta_m$ this reduces to

$$\sum_{j=0}^{m-1} \frac{H_n^j(\zeta_m)}{[j]_{\zeta_m}!} = \prod_{j=1}^{m-1} \left( I_n + c_j(\zeta_m) H_{n,j}(\zeta_m) \right) \tag{28}$$

because $H_n^m(\zeta_m) = H_{n,k}^{(m)} [m]_{\zeta_m}! = 0.$

Therefore we get from (24)

$$\left( \sum_{j=0}^{m-1} \frac{H_n^j(\zeta_m)}{[j]_{\zeta_m}!} \right)^{-1} P_n(\zeta_m) = \prod_{j\geq m} \left( I_n + c_j(\zeta_m) H_{n,j}(\zeta_m) \right). \tag{29}$$

Comparing with (26) we see that $c_j\left(\zeta_m\right)=0$ if $j \geq m$ is not a multiple of $m$ and that $c_{jm}\left(\zeta_m\right)=c_j$.

**Corollary 4.4** (Carlitz [2])

*Let $p$ be a prime number. If $n > p$ is relatively prime to $p$ then $c_n \equiv 0 \bmod p$. If $n = pm$ then $c_n = c_{pm} \equiv c_m \bmod p$.*

In another direction we prove a slight extension of the fact that $r_{2^n}(0)=1$ and $r_n(0)=0$ else.

**Theorem 4.5**

*The identity (16) reduces modulo $q^2$ to*

$$1+x+(1-q)x^2+(1-2q)x^3+\cdots \equiv \prod_{n\geq 1}\left(1+g_n(q)x^n\right) \bmod q^2 \tag{30}$$

*where $g_1(q)=1,\ \ g_{2^n}(q)=1-2^{n-1}q$ for $n>0$ and $g_{2n}(q)=0,\ \ g_{2n+1}(q)=-q$ else.*

**Proof**

For any commutative ring $R$ with identity the infinite product $\prod_{n\geq 1}\left(1+g_nx^n\right)$ with $g_n \in R$ can be expanded into a formal power series $1+a_1x+a_2x^2+\cdots$.

Let us choose $R=\mathbb{Z}[q]/\left(q^2\right)$. Its elements can be written as $a+bq$ with integers $a,b$ and $q^2=0$.

Since $[n]!=(1+q)\left(1+q+q^2\right)\cdots\left(1+q+\cdots+q^{n-1}\right)\equiv(1+q)^{n-1}\equiv 1+(n-1)q \bmod q^2$ we see that $\frac{1}{[n]!}\equiv 1-(n-1)q \bmod q^2$. Therefore $\exp_q(x)\equiv 1+x+(1-q)x^2+(1-2q)x^3+\cdots \bmod q^2$.

On the other hand we have

$$e_n(q)=\sum_{d|n,d>1}(-1)^d\frac{e_{n/d}(q)^d}{d}+\frac{(1-q)^{n-1}}{n[n]}.$$

Since

$$\frac{(1-q)^{n-1}}{n[n]}\equiv\frac{1-(n-1)q}{n(1+q)}\equiv\frac{\left(1-(n-1)q\right)(1-q)}{n}\equiv\frac{1-nq}{n}=\frac{1}{n}-q \bmod q^2$$

it suffices to show that $g_n(q)$ satisfies

$$\sum_{d|n}(-1)^d\frac{g_{n/d}^d(q)}{d}\equiv\left(q-\frac{1}{n}\right)\bmod q^2.$$

For $d=1$ we get $-g_n(q)$. For $d=n$ we get $\frac{(-1)^n}{n}$.

For $1<d<n$ we get $g_{n/d}^d=0$ except if $\frac{n}{d}=2^k$ for some $k$.

Thus for odd $n$ we get $-g_n(q)-\frac{1}{n}=q-\frac{1}{n}.$

Let now $n=2^k u$. Here we need only consider $d=1,\ d=n$ and $d=2^i u$ for $0\le i<k$.

We get $(-1)^d \frac{g_{n/d}^d(q)}{d}=(-1)^{2^i u}\frac{g_{2^{k-i}}^{2^i u}(q)}{2^i u}=(-1)^{2^i u}\frac{\left(1-2^{k-i-1}q\right)^{2^i u}}{2^i u}=(-1)^{2^i u}\left(\frac{1}{2^i u}-2^{k-i-1}q\right).$

This gives summed up

$$-\left(\frac{1}{u}-2^{k-1}q\right)+\left(\frac{1}{2u}-2^{k-2}q\right)+\left(\frac{1}{4u}-2^{k-3}q\right)+\cdots+\left(\frac{1}{2^{k-1}u}-q\right)=-\frac{1}{2^{k-1}u}+q$$

Together with $-g_n(q)+\frac{1}{n}=\frac{1}{2^k u}$ we get $q-\frac{1}{n}.$

Let us verify this for the first terms of $e_n(q)$:

$$e_1(q)=1=g_1(q),\ \ e_2(q)=\frac{1}{1+q}\equiv 1-q=1-2^0 q=g_2(q),$$

$$e_3(q)=-\frac{q}{1+q+q^2}\equiv -q\left(1-q\right)\equiv -q=g_3(q),$$

$$e_4(q)=\frac{1+q^2+q^3}{(1+q)(1+q+q^2+q^3)}\equiv(1-q)^2\equiv 1-2q=g_4(q),\cdots.$$